\documentclass[12pt]{amsart}
\usepackage[square,compress,comma, numbers,sort]{natbib}
\usepackage[colorlinks=true, citecolor=blue, linkcolor=blue]{hyperref}
\usepackage{amsfonts}
\allowdisplaybreaks[4]
\usepackage{graphicx}
\usepackage{epstopdf}

\usepackage{amssymb}
\usepackage{amsmath}
\usepackage{color}
\usepackage{enumerate}
\newcommand{\abs}[1]{\left\lvert #1 \right\rvert}

\def\E#1{\mathbb{E}\left \{#1 \right\}}

\definecolor{c20}{rgb}{0.,0.7,0.}
\definecolor{c30}{rgb}{0.,0.,1.}
\definecolor{c40}{rgb}{1,0.1,0.7}
\definecolor{c50}{rgb}{1,0,0}
\definecolor{c60}{rgb}{1,0.9,0.1}
\definecolor{c70}{rgb}{0.50,1.00,0.00}

\def\N{\mathbb{N}}
\numberwithin{equation}{section}
\newtheorem{theo}{Theorem}[section]
\newtheorem{sat}[theo]{Proposition}
\newtheorem{de}[theo]{Definition}
\newtheorem{lem}{Lemma}[section]

\newtheorem{korr}[theo]{Corollary}

\newtheorem{remarks}[theo]{Remarks}
\numberwithin{equation}{section}

\newcommand{\prooftheo}[1]{ \textsc{Proof of Theorem} \ref{#1} }

\newcommand{\prooflem}[1]{\textsc{Proof of Lemma} \ref{#1}}

\newcommand{\pk}[1]{\mathbb{P} \left\{ #1 \right\} }

\newcommand{\QED}{\hfill $\Box$}
\newcommand{\COM}[1]{}
\def\IF{\infty}
\newcommand{\R}{\mathbb{R}}

\newcommand{\BQN}{\begin{eqnarray}}
\newcommand{\EQN}{\end{eqnarray}}
\newcommand{\BQNY}{\begin{eqnarray*}}
\newcommand{\EQNY}{\end{eqnarray*}}
\def\ldot{, \ldots,}

\def\polhk#1{\setbox0=\hbox{#1}{\ooalign{\hidewidth
\lower1.5ex\hbox{`}\hidewidth\crcr\unhbox0}}}

\newcommand{\kb}[1]{\boldsymbol{#1}}
\newcommand{\vk}[1]{\kb{#1}}

\usepackage{mathtools}

\newcommand{\nelem}[1]{{Lemma \ref{#1}}}

\newcommand{\netheo}[1]{{Theorem \ref{#1}}}

\def\rw{\rightarrow}

\def\IF{\infty}

\date{}

\def\LT{\left}
\def\RT{\right}

\def\rw{\rightarrow}

\def\vn{\varepsilon}
\def\Var{\text{Var}}

\newcommand{\limit}[1]{\lim_{#1 \to \infty}}

\newcommand{\BS}{\begin{sat}}
\newcommand{\ES}{\end{sat}}
\newcommand{\BT}{\begin{theo}}
\newcommand{\ET}{\end{theo}}
\newcommand{\BK}{\begin{korr}}
\newcommand{\EK}{\end{korr}}

\newcommand{\BD}{\begin{de}}
\newcommand{\ED}{\end{de}}
\newcommand{\BIT}{\begin{itemize}}
\newcommand{\EIT}{\end{itemize}}
\newcommand{\BDI}{\begin{description}}
\newcommand{\EDI}{\end{description}}
\newcommand{\BRM}{\begin{remarks}}
\newcommand{\ERM}{\end{remarks}}
\newcommand{\BEL}{\begin{lem}}
\newcommand{\EEL}{\end{lem}}

\def\rw{\rightarrow}
\def\LT{\left}
\def\RT{\right}

\def\Var{\text{Var}}
\def\vn{\varepsilon}

\def\LT{\left}
\def\RT{\right}

\def\rw{\rightarrow}

\def\vn{\varepsilon}
\def\Var{\text{Var}}

\def\LT{\left}
\def\RT{\right}

\def\rw{\rightarrow}

\def\vn{\varepsilon}
\def\Var{\text{Var}}

\def\Z{\mathbb{Z}}

\def\v{\vk{v}}

\begin{document}
	
\title{Extremes of Locally-stationary Chi-square processes on discrete grids}

\author{Long Bai}
\address{Long Bai,
Department of Actuarial Science, University of Lausanne, UNIL-Dorigny, 1015 Lausanne, Switzerland
}
\email{Long.Bai@unil.ch}

	\bigskip
	
	\date{\today}
	\maketitle

{\bf Abstract}: For $X_i(t), i=1\ldot n, t\in [0,T]$ centered Gaussian processes, the chi-square process $\sum_{i=1}^{n}X_i^2(t)$ appears naturally as limiting processes in various statistical models. In this paper, we are concerned with the exact tail asymptotics of the supremum taken over discrete grids
of a class of locally stationary chi-square processes where $X_i(t),\ 1\leq i\leq n$ are not identical. An important tool for
establishing our results is a generalisation of Pickands lemma under the discrete scenario.  An application related to the change-point problem is discussed.

{{\bf Key Words:}  Chi-square processes; asymptotic methods; Pickands constant; change-point problem.
	

\section{Introduction and Main Result}

Numerous applications, especially from statistics are concerned with the supremum of chi-square processes over a discrete grid, which is threshold dependent.\\ 
\COM{We give more precise definition of the grids below: for given positive constants $\delta, \lambda $  we define $R_\lambda^\delta(u) = \delta u^{-\lambda } \cap \mathbb{Z}$ to be a grid indexed by the threshold $u$. Under the Pickands condition \eqref{lrr2}, the grid is said to be dense if $\lambda>\alpha   $ and  Pickands gird if $\alpha=\lambda$.}
The investigate of the extremes of chi-square process $\sum_{i=1}^nX_i^2(t)$
is initiated by the studies of high excursions of envelope of a Gaussian process, see e.g., \cite{BN1969} and generalized in \cite{Lindgren1980a,Lindgren1980b,Lindgren1989}. When $X_i's$ are stationary, \cite{Albin1990,Albin1992} develop the Berman's approach in \cite{Berman82} to obtain an asymptotic behavior of large deviation probabilities of the stationary chi-square processes. When $X_i, i=1\ldot n$ are locally-stationary Gaussian processes, \cite{LJ2017} obtains the extreme of the supremum of the locally-stationary Gaussian process. See \cite{ELPNT2017,EGCP2018} for more literature about locally stationary Gaussian processes.\\
In the case that the grid is very dense, then the tail asymptotics of supremum over $[0,T]$ or the dense grid is the same. The deference up to the constant appears for so-called Pickands grids, and in the case of rare grids, the asymptotics of supremum is completely different and usually very simple, see \cite{KP2018}.\\
Before giving the locally-stationary chi-square processes, we introduce  a class
of {\it locally stationary} Gaussian processes, considered by Berman in \cite{Berman92}, see also \cite{MR1342094,MR1062062,MR3679987, ChanLai,Marek,MR3638374}.
Specifically, let $X_i(t), i=1\ldot n, t\in [0,T]$ be centered Gaussian processes with unit variance and correlation function $r$ satisfying
 \BQN\label{locall-station}
\lim_{\epsilon\to0}\sup_{t,t+s\in[0,T], |s|<\epsilon} \abs{ \frac{1-r_i(t,t+s)} {\abs{s}^\alpha} - a_i(t) }=0,
\EQN
where $a_i(\cdot)$ are a continuous positive function on $[0,T]$ and $\alpha\in(0,2]$. \\
Then the locally-stationary chi-square process can be defined as
\BQNY
\chi^2(t)=\sum_{i=1}^{n}X_i^2(t),\ t\in[0,T].
\EQNY
\COM{In this paper, we investigate
\BQNY
\pk{\sup_{}}
\EQNY}
Let $\mathcal H_\alpha^\eta(a)$ be the Pickands constant defined for $\eta \ge 0$ and $a>0$ by
$$ \mathcal H_{\alpha}^\eta(a) = \limit{S} \frac 1{S} \E{  \sup_{ t\in   [0,S] \cap \eta \Z} e^{\sqrt{2a} B_\alpha(t)- a\abs{t}^\alpha}},$$
where $\eta \Z=\R$ if $\eta =0$. See \cite{PicandsA,Pit72, debicki2002ruin,DI2005,DE2014,DiekerY,DEJ14,Pit20, Tabis, DM, SBK, GeneralPit16} for various properties of $\mathcal{H}^\eta_{\alpha}(a)$.
\BT\label{theo-ls-tail}
Let $X_i(t),\ 1\leq i\leq n, \ t\in [0,T]$ be  centered, sample path continuous Gaussian processes with unit variance and correlation functions satisfying assumption \eqref{locall-station}. Suppose that $\eta_u,u>0$ are such that
\BQN \label{etau}
 \limit{u} \eta_u u^{1/\alpha} = \eta \in [0,\IF).
\EQN
If further $r(s,t)<1$ for all $s,t\in[0,T],s\neq t$, we have as $u\to\IF$
\BQN\label{concl-sta-tail}
\pk{ \sup_{t \in [0,T] \cap \eta_u  \Z}\chi^2(t)> u}  \sim  \int_0^{T}\int_{\v\in\mathbb{S}_{n-1}}\mathcal{H}^\eta_\alpha\LT(\sum_{i=1}^nv_i^2a_i(t) \RT) d\v dt
u^{1/\alpha} \pk{\chi^2(0)> u},
\EQN
where $\mathbb{S}_{n-1}=\{\v\in\R^{n}:\sum_{i=1}^n v_i^2=1\}$ .
\ET

Next in Section 2, we show an application related to change-point problem. The proof of \netheo{theo-ls-tail} and some lemmas are relegated to Section 3 and Section 4.

\section{Applications}

In \cite{ChCh2018}, they detect the change-point problem and give the generalized edge-count scan statistic expressed as
\BQNY
S(t)=X_1^2(t)+X_2^2(t),
\EQNY
where $X_i, i=1,2$ are two independent Gaussian processes which, respectively, have covariance functions
\BQNY
Cov\LT(X_1(s),X_1(t)\RT)=\frac{(s\wedge t)(1-(s\vee t))}{(s\vee t)(1-(s\wedge t))},\
Cov\LT(X_2(s),X_2(t)\RT)=\frac{(s\wedge t)(1-(s\vee t))}{\sqrt{(s\wedge t)(1-(s\wedge t))(s\vee t)(1-(s\vee t))}}.
\EQNY
Then in order to give the asymptotic $p$-value approximations, we need to investigate when $u$ is large enough for $0<T_1<T_2<\IF$
\BQNY
\pk{\sup_{t\in[T_1,T_2]\cap u^{-1}\Z}S(t)>u}.
\EQNY
\BT\label{application1}
we have as $u\rw\IF$
\BQNY
\pk{\sup_{t\in[T_1,T_2]\cap u^{-1}\Z}S(t)>u}\sim \frac{ue^{-u/2}}{2\pi}\int_{v_1^2+v_2^2=1}\int_{T_1}^{T_2}\mathcal{H}^1_1
\LT(\frac{v_1^2+\frac{1}{2}v_2^2}
{t(1-t)}\RT)dv_1d v_2dt.
\EQNY
\ET

\def\vv{\nu}
\def\LL{\rho}
\section{Proofs}
Below $\lfloor x \rfloor$ stands for the integer part of $x$
{and $\lceil x\rceil$ is the smallest integer not less than $x$}. Further $\Psi$ is the survival function of an $N(0,1)$ random variable.
During the following proofs, $\mathbb{Q}_i, i\in\N$ are some positive constants which can be different from line by line and for  interval $\Delta_1,\Delta_2  \subseteq[0,\IF)$ we denote
\BQNY
\mathcal{K}_u(\Delta_1):=\pk{\sup_{t\in\Delta_1}\chi^2(t)>u}, \ \
\mathcal{K}_u(\Delta_1,\Delta_2):=\pk{\sup_{t\in\Delta_1}\chi^2(t)>u,\sup_{t\in\Delta_2} \chi^2(t)>u}.
\EQNY

\def\Bale{\mathcal{B}_{\alpha,\lambda}^\eta}
\def\Seu{S_{\eta,u}}
\def\Se{S_\eta}

\prooftheo{theo-ls-tail}
For any $\theta>0$ and $\lambda>0$, set
\BQNY
&&I_k(\theta)=[k\theta,(k+1)\theta]\cap \eta_u \Z, \quad k\in \N,\quad N(\theta)=\LT\lfloor\frac{T}{\theta}\RT\rfloor,\\
&&J^k_l(u)=\LT[k\theta+lu^{-1/\alpha}\lambda,k\theta+(l+1)u^{-1/\alpha}\lambda\RT]\cap \eta_u \Z,\quad M(u)=\LT\lfloor\frac{\theta u^{1/\alpha}}{\lambda}\RT\rfloor,\\
&&K^k_l(u)=\LT[k\theta+lu^{-1/\alpha}\lambda,k\theta+(l+1)u^{-1/\alpha}\lambda\RT].
\EQNY
We have
\BQNY
\sum_{k=0}^{N(\theta)-1}\LT(\sum_{l=0}^{M(u)-1}\mathcal{K}_u(J^k_l(u))\RT)
-\sum_{i=1}^{4}\mathcal{A}_i(u)\leq\mathcal{K}_u([0,T]\cap \eta_u \Z)\leq\sum_{k=0}^{N(\theta)}
\mathcal{K}_u(I_k(\theta))
\leq\sum_{k=0}^{N(\theta)}\LT(\sum_{l=0}^{M(u)}\mathcal{K}_u(J^k_l(u))\RT),
\EQNY
where
\BQNY
&&\mathcal{A}_i(u)=\sum_{(k_1,l_1,k_2,l_2)\in \mathcal{L}_i}\mathcal{K}_u(K^{k_1}_{l_1}(u),K^{k_2}_{l_2}(u)),\ i=1,2,3,4,
\EQNY
with
\BQNY
&&\mathcal{L}_{1}=\LT\{0\leq k_1= k_2\leq N(\theta)-1,0\leq l_1+1=l_2\leq M(u)-1\RT\},\\
&&\mathcal{L}_{2}=\LT\{0\leq k_1+1= k_2\leq N(\theta)-1, l_1=M(u), l_2=0\RT\},\\
&&\mathcal{L}_3=\LT\{0\leq k_1+1<k_2\leq N(\theta)-1,0\leq l_1,l_2\leq M(u)-1\RT\},\\
&&\mathcal{L}_4=\LT\{0\leq k_1\leq k_2\leq N(\theta)-1,k_2-k_1\leq 1,0\leq l_1,l_2\leq M(u)-1\RT\}\setminus\LT(\mathcal{L}_1\cup \mathcal{L}_2\RT).
\EQNY
By \nelem{im1}, we have as $u\rw\IF,\ \lambda\rw\IF,\ \theta\rw 0$
\BQN\label{theoeq1}
\sum_{k=0}^{N(\theta)}\LT( \sum_{l=0}^{M(u)}\mathcal{K}_u(J^k_l(u))\RT)
&=&\sum_{k=0}^{N(\theta)}\LT(\sum_{l=0}^{M(u)}\pk{\sup_{t\in J^k_l(u)}\chi^2(t)>u}\RT) \nonumber\\
&\leq& \sum_{k=0}^{N(\theta)}\LT(\sum_{l=0}^{M(u)}
\int_{\v\in\mathbb{S}_{n-1}}\mathcal{H}^\eta_\alpha\LT(\sum_{i=1}^nv_i^2 (a_i(k\theta)+\vn_\theta) \RT) \lambda d\v \pk{\chi^2(0)>u}\RT)\nonumber\\
&\sim& \sum_{k=0}^{N(\theta)}\theta\int_{\v\in\mathbb{S}_{n-1}}\mathcal{H}^\eta_\alpha\LT(\sum_{i=1}^nv_i^2 (a_i(k\theta)+\vn_\theta) \RT) d\v u^{1/\alpha}\pk{\chi^2(0)>u} \nonumber\\
&\sim&\int_0^{T}\int_{\v\in\mathbb{S}_{n-1}}\mathcal{H}^\eta_\alpha\LT(\sum_{i=1}^nv_i^2a_i(t) \RT) d\v dt
u^{1/\alpha} \pk{\chi^2(0)> u}.
\EQN
Similarly, we have as $u\rw\IF,\ \lambda\rw\IF,\ \theta\rw 0$
\BQNY
\sum_{k=0}^{N(\theta)-1}\LT(\sum_{l=0}^{M(u)-1}\mathcal{K}_u(J^k_l(u))\RT)
\geq \int_0^{T}\int_{\v\in\mathbb{S}_{n-1}}\mathcal{H}^\eta_\alpha\LT(\sum_{i=1}^nv_i^2a_i(t) \RT) d\v dt
u^{1/\alpha} \pk{\chi^2(0)> u}.
\EQNY
Next we focus on the analysis of $\mathcal{A}_i(u), i=1,2$. For $(k_1,l_1,k_2,l_2)\in \mathcal{L}_{1}$, without loss of generality, we assume $k_1=k_2$ and  $l_1+1 = l_1$. Then set
\BQNY
&&\LT(K^{k_1}_{l_1}(u)\RT)^1=\LT[k_1\theta+l_1u^{-1/\alpha}\lambda,
k_1\theta+(l_1+1)u^{-1/\alpha}(\lambda-\sqrt{\lambda})\RT],\\
&&\LT(K^{k_1}_{l_1}(u)\RT)^2=\LT[k_1\theta+(l_1+1)u^{-1/\alpha}(\lambda-\sqrt{\lambda}),
k\theta+(l+1)u^{-1/\alpha}\lambda\RT].
\EQNY
Then we have
\BQNY
\mathcal{A}_1(u)\leq\sum_{(k_1,l_1,k_2,l_2)\in \mathcal{L}_i}\LT(\mathcal{K}_u\LT(\LT(K^{k_1}_{l_1}(u)\RT)^{1},K^{k_2}_{l_2}(u)\RT)
+ \mathcal{K}_u\LT(\LT(K^{k_1}_{l_1}(u)\RT)^{2}\RT)\RT).
\EQNY
Analogously as in \eqref{theoeq1}, we have as $u\rw\IF,\ \lambda\rw\IF,\ \theta\rw 0$
\BQNY
\sum_{(k_1,l_1,k_2,l_2)\in \mathcal{L}_1}\mathcal{K}_u\LT(\LT(K^{k_1}_{l_1}(u)\RT)^{2}\RT)
&\leq& \sum_{k_1=0}^{N(\theta)-1}\sum_{l_1=0}^{M(u)-1}
\mathcal{K}_u\LT(\LT(K^{k_1}_{l_1}(u)\RT)^{2}\RT)\\
&\leq& \sum_{k=0}^{N(\theta)}\LT(\sum_{l=0}^{M(u)}
\int_{\v\in\mathbb{S}_{n-1}}\mathcal{H}^\eta_\alpha\LT(\sum_{i=1}^nv_i^2 (a_i(k\theta)+\vn_\theta) \RT) \sqrt{\lambda} d\v \pk{\chi^2(0)>u}\RT)\nonumber\\
&\sim&\frac{1}{\sqrt{\lambda}}\int_0^{T}\int_{\v\in\mathbb{S}_{n-1}}\mathcal{H}^\eta_\alpha\LT(\sum_{i=1}^nv_i^2a_i(t) \RT) d\v dt
u^{1/\alpha} \pk{\chi^2(0)> u},\\
&=&o\LT(u^{1/\alpha} \pk{\chi^2(0)> u}\RT).
\EQNY
\COM{Then by \nelem{}, we have
\BQNY
\sum_{(k_1,l_1,k_2,l_2)\in \mathcal{L}_i}\mathcal{K}_u\LT(\LT(K^{k_1}_{l_1}(u)\RT)^{1},K^{k_2}_{l_2}(u)\RT)
\leq \mathbb{Q} \lambda^2
\EQNY}
Further, by \nelem{im1}
\BQNY
\mathcal{A}_1(u)
&\leq&\sum_{k=0}^{N(\theta)-1}\LT(\sum_{l=0}^{M(u)-1}\left(\mathcal{K}_u(K^k_l(u))
+\mathcal{K}_u(K^k_{l+1}(u))-\mathcal{K}_u(K^k_l(u)\cup K^k_{l+1}(u))\right)\RT)\nonumber\\
&\sim &  \sum_{k=0}^{N(\theta)-1}\LT(\LT(\mathcal{H}_\alpha[0,(a_k+\vn_\theta)^{\frac{1}{\alpha}}d^{-1/\alpha}\lambda]
+\mathcal{H}_\alpha[0,(a_k+\vn_\theta)^{\frac{1}{\alpha}}d^{-1/\alpha}\lambda]
-\mathcal{H}_\alpha[0,2(a_k-\vn_\theta)^{\frac{1}{\alpha}}d^{-1/\alpha}\lambda]\RT)\RT.\\
&&\LT.\times\sum_{l=0}^{M(u)-1}\pk{\chi^2(0)> u}\RT)\nonumber\\
&\leq&\mathbb{Q}_1\LT(\sum_{k=0}^{N(\theta)-1}\LT((a_k+\vn_\theta)^{\frac{1}{\alpha}}-(a_k
-\vn_\theta)^{\frac{1}{\alpha}}\RT)\theta\RT)u^{2/\alpha^*}\pk{\chi^2(0)> u}\\
&=&o\left(u^{2/\alpha^*}\pk{\chi^2(0)> u}\right), \ u\rw\IF,\ \ \lambda\rw\IF,\theta\rw0.
\EQNY
Similarly, by \nelem{im1}
\BQNY
\mathcal{A}_2(u)
&=&\sum_{k=0}^{N(\theta)-1}\mathcal{K}_u(J^{k}_{M(u)-1}(u), J^{k+1}_{0}(u))\\
&\leq&\sum_{k=0}^{N(\theta)-1}\pk{\sup_{t\in [0,2\lambda] }\chi^2((k+1)\theta-u^{-2/\alpha}t)>u,\sup_{t\in [0,2\lambda]}\chi^2((k+1)\theta+u^{-2/\alpha}t)>u}\\
&=&\sum_{k=0}^{N(\theta)-1}\LT(\pk{\sup_{t\in [0,2\lambda] }\chi^2((k+1)\theta-u^{-2/\alpha}t)>u}+\pk{\sup_{t\in [0,2\lambda]}\chi^2((k+1)\theta+u^{-2/\alpha}t)>u}\RT.\\
&&\LT.-\pk{\sup_{t\in [-2\lambda,2\lambda] }\chi^2((k+1)\theta-u^{-2/\alpha}t)>u}\RT)\nonumber\\
&\sim &  \sum_{k=0}^{N(\theta)-1}\LT(\LT(2\mathcal{H}_\alpha[0,2(a_{k+1}+\vn_\theta)^{\frac{1}{\alpha}}d^{-1/\alpha}\lambda]
-\mathcal{H}_\alpha[-2(a_k-\vn_\theta)^{\frac{1}{\alpha}}d^{-1/\alpha}\lambda,
2(a_k-\vn_\theta)^{\frac{1}{\alpha}}d^{-1/\alpha}\lambda]\RT)\RT.\\
&&\LT.\times\sum_{l=0}^{M(u)-1}\pk{\chi^2(0)>u}\RT)\nonumber\\
&\leq&\mathbb{Q}_2\LT(\sum_{k=0}^{N(\theta)-1}\LT((a_k+\vn_\theta)^{\frac{1}{\alpha}}-(a_k
-\vn_\theta)^{\frac{1}{\alpha}}\RT)\theta\RT)u^{2/\alpha}
\pk{\chi^2(0)>u}\\
&=&o\left(u^{2/\alpha}\pk{\chi^2(0)>u}\right), \ u\rw\IF,\ \ \lambda\rw\IF,\theta\rw0.
\EQNY

For any $\theta>0$
\BQNY
\E{X_i(t)X_i(s)}=r(s,t)\leq 1-\delta(\theta)
\EQNY
for $(s,t)\in J^{k_1}_{l_1}(u)\times J^{k_2}_{l_2}(u),(j_1,k_1,j_2,k_2)\in\mathcal{L}_3$ where $\delta(\theta)>0$ is related to $\theta$. Then we have
\BQNY
\mathcal{A}_3(u)
&\leq&N(\theta)M(u)2\Psi\LT(\frac{2u-\mathbb{Q}_3}{d\sqrt{4-\delta(\theta)}}\RT)\\
&\leq&\frac{T}{\lambda}u^{2/\alpha}2\Psi\LT(\frac{2u-\mathbb{Q}_3}{d\sqrt{4-\delta(\theta)}}\RT)\\
&=&o\LT(u^{2/\alpha}\pk{\chi^2(0)>u}\RT), \ u\rw\IF, \lambda\rw\IF, \theta\rw 0.
\EQNY
where $\mathbb{Q}_3$ is a large constant.
Finally by \nelem{in1} for $u$ large enough and $\theta$ small enough
\BQNY
\mathcal{A}_4(u)
&\leq &\sum_{k=0}^{N(\theta)-1}\LT(\sum_{l=0}^{2M(u)}\sum_{i=2}^{2M(u)}
\mathcal{K}_u(J^k_l(u),J^{k}_{l+i}(u))\RT)\\
&\leq&\sum_{k=0}^{N(\theta)-1} \sum_{l=0}^{2M(u)}\pk{\chi^2(0)>u}
\LT(\sum_{i=1}^{\IF}\mathbb{Q}_4\exp\LT(-\frac{\mathbb{Q}_5}{8}\abs{i\lambda}^\alpha\RT)\RT)\\
&\leq&\mathbb{Q}_6\frac{T}{\lambda}u^{-2/\alpha}\pk{\chi^2(0)>u}\LT(\sum_{i=1}^{\IF}\exp\LT(-\frac{\mathbb{Q}_5}{8}\abs{i\lambda}^\alpha\RT)\RT)\\
&=&o\LT(u^{2/\alpha}\pk{\chi^2(0)>u}\RT), \ u\rw\IF, \lambda\rw\IF, \theta\rw 0.
\EQNY
Thus the claim follows.
\QED

\section{Appendix} 
\prooftheo{application1}  
For the correlation function $r_i(t)$ of $X_i(t),i=1,2$, a simple calculation shows that
\BQNY
\lim_{\epsilon\to0}\sup_{t,t+s\in[T_1,T_2], |s|<\epsilon} \abs{ \frac{1-r_i(t,t+s)} {\abs{s}} - a_i(t)}=0, \ i=1,2,
\EQNY
where $a_1(t)=\frac{1}{t(1-t)}$ and $a_2(t)=\frac{1}{2t(1-t)}$.\\
Then by \netheo{theo-ls-tail}, the result fellows.
\QED

\BEL\label{im1}
Let $X_i(t),\ 1\leq i\leq n, \ t\in [0,T]$ be  centered, sample path continuous Gaussian processes with unit variance and correlation functions satisfying assumption \eqref{locall-station}. Suppose that $\eta_u,u>0$ are such that
\BQN \label{etau}
 \limit{u} \eta_u u^{1/\alpha} = \eta \in [0,\IF).
\EQN
Set $ a:=a(t_0), \ t_0\in\R,$ and $K_u$ a family of index sets. If $\lim_{u\rw\IF}\sup_{k\in K_u}\abs{k u^{-2/\alpha}}\leq \theta$ for some small enough $\theta\geq0$,
we have for some constant $S>0$, $S_1,S_2\in\R$ with $S_1<S_2$ when $u$ large enough
\BQNY
\int_{\v\in\mathbb{S}_{n-1}}\mathcal{H}^\eta_\alpha\LT(\sum_{i=1}^nv_i^2(a_i-\vn_\theta) \RT)d\v &\leq&\lim_{u\rw\IF}\forall_{k\in K_u}\frac{\pk{ \sup_{t \in [S_1,S_2] \cap \eta_u  \Z}\chi^2(t_0+u^{-1/\alpha}kS+t)> u}}{\pk{\chi^2(0)> u}}\\
&\leq&  \int_{\v\in\mathbb{S}_{n-1}}\mathcal{H}^\eta_\alpha\LT(\sum_{i=1}^nv_i^2(a_i+\vn_\theta) \RT)d\v.
\EQNY
\EEL
\prooflem{im1}
 For $u$ large enough,  using the duality property of norm we find
 \BQNY
 &&\pk{ \sup_{t \in [S_1,S_2] \cap \eta_u  \Z}\chi^2(t_0+u^{-1/\alpha}kS+u^{-1/\alpha}t)> u}\\
 &&=
 \pk{ \sup_{(t,\v) \in ([S_1,S_2] \cap \eta_u  \Z)\times(\mathbb{S}^{d-1})}Y_u(t,\v)> u^{1/2}},
 \EQNY
 where $Y_u(t,\v)=\sum_{i=1}^{n}v_iX_i(t_0+u^{-1/\alpha}kS+u^{-1/\alpha}t)$.
 We know that the variance function of $Y_u(t,\v)$ always equal to 1 over $([S_1,S_2] \cap \eta_u  \Z)\times(\mathbb{S}^{d-1})$.
 The fields $Y_u(t,\v)$ can be represented as
\BQNY
Y_u(t,\widetilde{\vk{v}})=\sum_{i=2}^{n}v_i X_i(t_0+u^{-1/\alpha}kS+u^{-1/\alpha}t)
+\LT(1-\sum_{i=2}^{n}v_i^2\RT)^{1/2}  X_1(t_0+u^{-1/\alpha}kS+u^{-1/\alpha}t),
 \widetilde{\vk{v}}=(v_2,\cdots,v_{n}),
\EQNY
which is defined in $[-S_1,S_2]\times\widetilde{\mathbb{S}}^{d-1}$ where
$$\widetilde{\mathbb{S}}^{d-1}
=\LT\{\widetilde{\vk{v}}:\LT(\LT(1-\sum_{i=2}^{n}v_i^q \RT)^{1/q},v_2,\cdots,v_{n},\RT)\in\mathbb{S}^{d-1}\RT\}.$$
Furthermore, following the arguments as in \cite{Pitchi1994} we conclude that
 the correlation function $r_u(t,\widetilde{\vk{v}},s,\widetilde{\vk{w}})$ of $Y_u(t, \widetilde{\vk{v}})$ satisfies for $u$ large enough
\BQNY
&& r_u(t,\widetilde{\vk{v}},s,\widetilde{\vk{w}})\geq 1- u^{-2}\LT(\sum_{i=1}^nv_i^2(a_i+\vn_\theta)\RT)\abs{t-s}^\alpha
-\frac{1+\vn_\theta}{2}\sum_{i=2}^n(v_i-w_i)^2,\\
&& r_u(t,\widetilde{\vk{v}},s,\widetilde{\vk{w}})\leq 1- u^{-2}\LT(\sum_{i=1}^nv_i^2(a_i-\vn_\theta)\RT)\abs{t-s}^\alpha
-\frac{1-\vn_\theta}{2}\sum_{i=2}^n(v_i-w_i)^2.
\EQNY
Then the proof follows by similar arguments as in the proof of \cite{PL2015} [Theorem 6.1] with the case $\mu=\nu$. Consequently, we get
\BQNY
\pk{ \sup_{(t,\v) \in ([S_1,S_2] \cap \eta_u  \Z)\times(\mathbb{S}^{d-1})}Y_u(t,\v)> u^{1/2}}
&=&\pk{\sup_{(t,\widetilde{\v}) \in ([S_1,S_2] \cap \eta_u  \Z)\times(\widetilde{\mathbb{S}}^{d-1})}Y_u(t,\widetilde{\v})> u^{1/2}}\\
&\leq&\int_{\v\in\mathbb{S}_{n-1}}\mathcal{H}^\eta_\alpha\LT(\sum_{i=1}^nv_i^2(a_i+\vn_\theta) \RT)d\v \pk{\chi^2(0)> u},
\EQNY
and
\BQNY
\pk{ \sup_{(t,\v) \in ([S_1,S_2] \cap \eta_u  \Z)\times(\mathbb{S}^{d-1})}Y_u(t,\v)> u^{1/2}}
&\geq&\int_{\v\in\mathbb{S}_{n-1}}\mathcal{H}^\eta_\alpha\LT(\sum_{i=1}^nv_i^2(a_i-\vn_\theta) \RT)d\v \pk{\chi^2(0)> u}.
\EQNY
\QED

\BEL\label{in1}
Assume the same assumptions as in \nelem{im1}.  Further, let $\vn_0$ be such that for all $s,t\in[t_0-\vn_0,t_0+\vn_0]$ and $1\leq i\leq n$
\BQNY
\frac{a}{2}\abs{t-s}^{\alpha}\leq1-r_i(s,t)\leq2a\abs{t-s}^\alpha.
\EQNY
Then we can find a constant $\mathbb{C}$ such that for all $S>0$ and $ T_2-T_1>S$,
\BQNY
{\limsup}_{u\rw\IF}\sup_{k\in K_u}\frac{\pk{\mathcal{A}_1(u_k),\mathcal{A}_2(u_k)}}
{\pk{\chi^2(t_0)>u}}\leq \mathbb{C}\exp\LT(-\frac{a}{8}|T_2-T_1-S|^\alpha\RT),
\EQNY
where $\mathcal{A}_i(u_k)=\{\sup_{t\in [T_i,T_i+S]}\chi^2(u^{-2/\alpha}(t+kS)+t_0)>u\},\ i=1,2$, and
\BQNY
\lim_{u\rw\IF}\sup_{k\in K_u}\abs{u^{-2/\alpha}k S}\leq \vn_0.
\EQNY
\EEL
\prooflem{in1}
Through this proof, $\mathbb{C}_i, i\in\N$ are some positive constant.\\
Set $Y_u(t,\vk{v})=\sum_{i=1}^n v_iX_{i}(u^{-2/\alpha}(t+kS)+t_0),(t,\vk{v})\in\R\times\mathcal{S}_q$ which is a centered Gaussian field and $\mathcal{S}_q^{\delta}=\{\vk{v}\in\mathcal{S}_q:1-\sum_{i=1}^nv_i^2\leq \delta\},\delta>0$.\\
Below for $\Delta_1,\Delta_2\subseteq \R^{n+1}$, denote
\BQNY
\mathcal{Y}_u(\Delta_1,\Delta_2)=\pk{\sup_{(t,\vk{v})\in\Delta_1}Y_u(t,\vk{v})>u_k
,\sup_{(t,\vk{v})\in\Delta_2}Y_u(t,\vk{v})>u_k}.
\EQNY
We have
\BQNY
\pk{\mathcal{A}_1(u_k),\mathcal{A}_2(u_k)}&\geq&
\mathcal{Y}_u([T_1,T_1+S]\times\mathcal{S}^\delta_q,[T_2,T_2+S]\times\mathcal{S}^\delta_q),\\
\pk{\mathcal{A}_1(u_k),\mathcal{A}_2(u_k)}&\leq&
\mathcal{Y}_u([T_1,T_1+S]\times\mathcal{S}^\delta_q,[T_2,T_2+S]\times\mathcal{S}^\delta_q)\\
&&+
\mathcal{Y}_u([T_1,T_1+S]\times\mathcal{S}^\delta_q,[T_2,T_2+S]\times(\mathcal{S}_q
\setminus\mathcal{S}^\delta_q))\\
&&+\mathcal{Y}_u([T_1,T_1+S]\times(\mathcal{S}_q\setminus\mathcal{S}^\delta_q)
,[T_2,T_2+S]\times\mathcal{S}^\delta_q),
\EQNY
and
\BQNY
\mathcal{Y}_u([T_1,T_1+S]\times\mathcal{S}^\delta_q,[T_2,T_2+S]\times(\mathcal{S}_q
\setminus\mathcal{S}^\delta_q))
&\leq& \pk{
\sup_{(t,\vk{v})\in[T_2,T_2+S]\times(\mathcal{S}_q\setminus\mathcal{S}^\delta_q)}Y_u(t,\vk{v})>u_k}\\
&\leq& \exp\LT(-\frac{\LT(u_k-\mathbb{C}_1\RT)^2}{2(d^2-\delta)}\RT)\\
&=&o\LT(\pk{Z(t_0)>u_k}\RT),
\EQNY
as $u\rw\IF$ where the last second inequality follows from Borell inequality and the fact that
\BQNY
\sup_{(t,\vk{v})\in[T_2,T_2+S]\times(\mathcal{S}_q\setminus\mathcal{S}^\delta_q)}
\Var(Y_u(t,\vk{v}))=\sup_{\vk{v}\in(\mathcal{S}_q\setminus\mathcal{S}^\delta_q)}
\LT(\sum_{i=1}^n d_i^2v_i^2\RT)\leq d^2-\delta.
\EQNY
Similarly, we have
\BQNY
\mathcal{Y}_u([T_1,T_1+S]\times(\mathcal{S}_q\setminus\mathcal{S}^\delta_q)
,[T_2,T_2+S]\times\mathcal{S}^\delta_q)=o\LT(\pk{Z(t_0)>u_k}\RT),\ u\rw\IF.
\EQNY
Then we just need to focus on
$$\Pi(u):=\mathcal{Y}_u([T_1,T_1+S]\times\mathcal{S}^\delta_q,[T_2,T_2+S]\times\mathcal{S}^\delta_q).$$
We split $\mathcal{S}^\delta_q$ into sets of small diameters $\{\partial \mathcal{S}_i, 0\leq i\leq \mathcal{N}^*\}$, where
\BQNY
\mathcal{N}^*=\sharp\{\partial\mathcal{S}_i\}<\IF.
\EQNY
Further, we see that $\Pi(u)\leq \Pi_1(u)+\Pi_2(u)$ with
\BQNY
\Pi_1(u)=\underset{\partial\mathcal{S}_i\cap\partial\mathcal{S}_l=\emptyset}
{\sum_{0\leq i,l\leq \mathcal{N}^*}}
\mathcal{Y}_u([T_1,T_1+S]\times\partial\mathcal{S}_i,[T_2,T_2+S]\times\partial\mathcal{S}_l),\\
\Pi_2(u)=\underset{\partial\mathcal{S}_i\cap\partial\mathcal{S}_l\neq\emptyset}{\sum_{0\leq i,l\leq \mathcal{N}^*}}
\mathcal{Y}_u([T_1,T_1+S]\times\partial\mathcal{S}_i,[T_2,T_2+S]\times\partial\mathcal{S}_l),
\EQNY
where ${\partial\mathcal{S}_i\cap\partial\mathcal{S}_l\neq\emptyset}$ means $\partial\mathcal{S}_i,\partial\mathcal{S}_l$ are identical or adjacent, and ${\partial\mathcal{S}_i\cap\partial\mathcal{S}_l=\emptyset}$ means $\partial\mathcal{S}_i,\partial\mathcal{S}_l$ are neither identical nor adjacent.
Denote the distance of two set $\vk{A}, \vk{B}\in\R^n$ as
\BQNY
\rho(\vk{A},\vk{B})=\inf_{\vk{x}\in \vk{A},\vk{y}\in \vk{B}}\|\vk{x}-\vk{y}\|_2.
\EQNY
if $\partial\mathcal{S}_i\cap\partial\mathcal{S}_l=\emptyset$, then there exists some small positive constant $\rho_0$ (independent of $i,l$) such that $\rho(\partial\mathcal{S}_i,\partial\mathcal{S}_l)>\rho_0$.
Next we estimate $\Pi_1(u)$. For any $u\geq 0$
\BQNY
\Pi_1(u)\leq \pk{\underset{\vk{v}\in\partial\mathcal{S}_i,\vk{w}\in\partial\mathcal{S}_i}
{\sup_{(t,s)\in[T_1,T_1+S]\times[T_2,T_2+S]}}Z_u(t,\vk{v},s,\vk{w})>2u_k},
\EQNY
where $Z_u(t,\vk{v},s,\vk{w})=Y_u(t,\vk{v})+Y_u(s,\vk{w}), \ t,s\geq 0, \vk{v},\vk{w}\in\R^n$.\\
When  $u$ is sufficiently large for $(t,s)\in[T_1,T_1+S]\times[T_2,T_2+S],\vk{v}\in\partial\mathcal{S}_i\subset[-2,2]^n,
\vk{w}\in\partial\mathcal{S}_i\subset[-2,2]^n$, with $\rho(\partial\mathcal{S}_i,\partial\mathcal{S}_l)>\rho_0$ we have
\BQNY
Var(Z_u(t,\vk{v},s,\vk{w}))&\leq&\sum_{i=1}^n(v_i^2+w_i^2+2v_iw_i)d_i^2\\
&\leq&4d^2-2\sum_{i=1}^n(v_i-w_i)^2d_i^2\\
&=&4d^2-2d_n^2\rho_0\\
&\leq&d^2(4-\delta_0),
\EQNY
for some $\delta_0>0$. Therefore, it follows from the Borell inequality that
\BQNY
\Pi_1(u)\leq \mathbb{C}_2 \mathcal{N}^*\exp\LT(-\frac{\LT(2u_k-\mathbb{C}_3\RT)^2}{2d^2(4-\delta_0)}\RT)
&=&o\LT(\pk{Z(t_0)>u_k}\RT),
\ u\rw\IF,
\EQNY
with
$$\mathbb{C}_3=\E{\underset{(\vk{v},\vk{w})\in[-2,2]^2n}{\sup_{(t,s)\in[T_1,T_1+S]\times[T_2,T_2+S]}}Z_u(t,\vk{v},s,\vk{w})}<\IF.$$
Now we consider $\Pi_2(u)$. Similar to the argumentation as in {\bf Step1} of the proof of \nelem{im1}. we set $\widetilde{Y}_u(t,\widetilde{\vk{v}})=Y_u(t,Q\widetilde{\vk{v}})$ and $\widetilde{Z}_u(t,\widetilde{\vk{v}},s,\widetilde{\vk{w}})=\widetilde{Y}_u(t,\widetilde{\vk{v}})+\widetilde{Y}_u(s,\widetilde{\vk{w}})$ with $\widetilde{\vk{v}},\widetilde{\vk{w}}\in \R^{n-1}$.
Since for $(t,s)\in[T_1,T_1+S]\times[T_2,T_2+S],\widetilde{\vk{v}}\in[-2,2]^{n-1},
\widetilde{\vk{w}}\in[-2,2]^{n-1}$, we have
\BQNY
2d^2\leq Var(\widetilde{Z}_u(t,\widetilde{\vk{v}},s,\widetilde{\vk{w}}))&\leq&\sum_{i=1}^n(v_i^2+w_i^2+2r(u^{-2/\alpha}(t+kS)+t_0,
u^{-2/\alpha}(s+kS)+t_0)v_i w_i)d_i^2\\
&\leq&2d^2+2\LT(1-\frac{a}{2}u^{-2}\abs{t-s}^\alpha\RT)\sum_{i=1}^n v_iw_id_i^2\\
&\leq&4d^2-d^2au^{-2}\abs{t-s}^\alpha\\
&\leq&4d^2-d^2au^{-2}\abs{T_2-T_1-S}^\alpha.
\EQNY
Set
\BQNY
\overline{Z}_u(t,\widetilde{\vk{v}},s,\widetilde{\vk{w}})=
\frac{\widetilde{Z}_u(t,\widetilde{\vk{v}},s,\widetilde{\vk{w}})}
{Var(Z_u(t,\widetilde{\vk{v}},s,\widetilde{\vk{w}}))}.
\EQNY
Borrowing the arguments of the proof in \cite{Pit96} [Lemma 6.3] we show that
\BQNY
\E{\LT(\overline{Z}_u(t,\widetilde{\vk{v}},s,\widetilde{\vk{w}})-\overline{Z}_u(t',\widetilde{\vk{v}'},s',\widetilde{\vk{w}'})\RT)}
\leq 4\LT(\E{(\widetilde{Y}_u(t,\widetilde{\vk{v}})-\widetilde{Y}_u(t',\widetilde{\vk{v}'}))^2}
+\E{(Y_u(s,\widetilde{\vk{w}})-Y_u(s',\widetilde{\vk{w}'}))^2}\RT).
\EQNY
Moreover,
\BQNY
r_1(t,\widetilde{\vk{v}},s,\widetilde{\vk{v}}')=1- u^{-2}a(t-s)^\alpha
-\frac{1}{2}\sum_{i=2}^nd_i^2(v_i-v_i')^2+o\LT(\sum_{i=2}^nd_i^2(v_i-v_i')^2+u^{-2}\RT), \widetilde{\vk{v}},\widetilde{\vk{v}}'\rw\widetilde{\vk{0}}, u\rw\IF.
\EQNY

Then we have
\BQNY
\E{(Y_u(t,\widetilde{\vk{v}})-Y_u(t',\widetilde{\vk{v}'}))^2}\leq 4d^2a u^{-2}\abs{t-t'}^\alpha+2\sum_{i=2}^n(v_i-v_i')^2.
\EQNY
Therefore
\BQNY
&&\E{\LT(\overline{Z}_u(t,\widetilde{\vk{v}},s,\widetilde{\vk{w}})
-\overline{Z}_u(t',\widetilde{\vk{v}'},s',\widetilde{\vk{w}'})\RT)}\\
&&\leq 16d^2a u^{-2}\abs{t-t'}^\alpha+16d^2a u^{-2}\abs{s-s'}^\alpha+8\sum_{i=2}^n(v_i-v_i')^2
+8\sum_{i=2}^n(w_i-w_i')^2.
\EQNY
Set $\zeta(t,s,\widetilde{\vk{v}},\widetilde{\vk{w}}), t,s\geq0, \vk{v},\vk{w}\in\R^{n-1}$ is a stationary Gaussian field with unit variance and correlation function
\BQNY
r_{\zeta}(t,s,\widetilde{\vk{v}},\widetilde{\vk{w}})=\exp\LT(-9 d^2at^\alpha-9 d^2as^\alpha-5\sum_{i=2}^nv_i^2-5\sum_{i=2}^nw_i^2\RT).
\EQNY
Then
\BQNY
\Pi_2(u)&\leq& \pk{\underset{\widetilde{\vk{v}}\in Q^{-1}\mathcal{S}_q,\widetilde{\vk{w}}\in Q^{-1}\mathcal{S}_q}
{\sup_{(t,s)\in[T_1,T_1+S]\times[T_2,T_2+S]}}\widetilde{Z}_u(t,\widetilde{\vk{v}},s,\widetilde{\vk{w}})>2u_k}\\
&\leq&\pk{\underset{\widetilde{\vk{v}}\in Q^{-1}\mathcal{S}_q,\widetilde{\vk{w}}\in Q^{-1}\mathcal{S}_q}
{\sup_{(t,s)\in[T_1,T_1+S]\times[T_2,T_2+S]}}\zeta(u^{-2/\alpha}t,u^{-2/\alpha}s,\widetilde{\vk{v}},\widetilde{\vk{w}})
>\frac{2u_k}{\sqrt{4d^2-d^2au^{-2}\abs{T_2-T_1-S}^\alpha}}}.
\EQNY
Then following the similar argumentation as in \cite{EnkelejdJi2014Chi}, we have
\BQNY
\Pi_2(u)\leq \mathbb{C}_4u_k^{M-2}\exp\LT(-\frac{u_k^{2}}{2d^2}-\frac{a}{8}\abs{T_2-T_1-S}^\alpha\RT)
\EQNY
where $M=0$ when $p\in(1,2)\cup(2,\IF]$ and $M=m$ when $p=2$.
Thus we have
\BQNY
\limsup_{u\rw\IF}\frac{\Pi_2(u)}{\pk{Z(t_0)>u_k}}\leq \mathbb{C}_5\exp\LT(-\frac{a}{8}|T_2-T_1-S|^\alpha\RT).
\EQNY
Thus we complete the proof.
\COM{
Set $Y(t,\vk{v})=\sum_{i=1}^n d_iv_iX_{i}(t),(t,\vk{v})\in\R\times\mathcal{S}_q$ which is a centered Gaussian field.
First we have for $u$ large
\BQNY
&&2\leq4-4a|t-s|^\alpha\leq\Var(X_i(s)+X_i(t))=2+2r(|t-s|)\leq4-a|t-s|^\alpha
\leq4-a\abs{T_2-T_1-S}\alpha u^{-2/c},\\
&&1-Cov(X_i(s)+X_i(t),X_i(s')+X_i(t'))\leq 1-r(|s-s'|)+1-r(|t-t'|)\leq 2a|s-s'|^{\alpha}+2a|t-t'|^{\alpha},
\EQNY
with $s\in[t_0+T_1u^{-2/\alpha^*},t_0+(T_1+S)u^{-2/\alpha^*}],
t\in[t_0+T_2u^{-2/\alpha^*},t_0+(T_2+S)u^{-2/\alpha^*}]$.\\
Set $\mathcal{S}_q=\{\vk{v}\in\R^n:||\vk{v}||_q=1\}, \mathcal{S}'_q=\{\vk{v}'\in\R^n:||\vk{v}'||_q=1\}$. We divide $\mathcal{S}_q$ and $\mathcal{S}_q$ respectively into $2^n$ parts according to $v_i\geq 0$ and $v_i\leq 0$  which are denote by $\mathcal{S}_q(j),\mathcal{S}'_q(j), j=1,\cdots,2^n$.
We know
\BQNY
\pk{\mathcal{A}_1(u_k),\mathcal{A}_2(u_k)}
&=&\pk{\sup_{(t,\vk{v})\in[T_1,T_1+S]\times\mathcal{S}_q}Y(u^{-2/\alpha^*}t+t_0,\vk{v})>u^{1/c}_k
,\sup_{(s,\vk{v}')\in[T_2,T_2+S]\times\mathcal{S}'_q}Y(u^{-2/\alpha^*}s+t_0,\vk{v}')>u^{1/c}_k}\\
&\leq&\pk{\sup_{(t,s,\vk{v},\vk{v}')\in[T_1,T_1+S]\times[T_2,T_2+S]\times\mathcal{S}_q
\times\mathcal{S}'_q}Y(u^{-2/\alpha^*}t+t_0,\vk{v})+Y(u^{-2/\alpha^*}s+t_0,\vk{v}')>2u^{1/c}_k}\\
&\leq&\sum_{j=1}^{2^n}\sum_{k=1}^{2^n}\pk{\sup_{(t,s,\vk{v},\vk{v}')\in[T_1,T_1+S]\times[T_2,T_2+S]\times\mathcal{S}_q(j)
\times\mathcal{S}'_q(k)}Y(u^{-2/\alpha^*}t+t_0,\vk{v})+Y(u^{-2/\alpha^*}s+t_0,\vk{v}')>2u^{1/c}_k}.
\EQNY
Next we consider one term in the sum of the upper inequality.
Without lose of generality, we set $\mathcal{S}_q(1)=\{\vk{v}\in\R^n:||\vk{v}||_q=1, v_i\geq 0\}$ and $\mathcal{S}'_q(1)=\{\vk{v}'\in\R^n:||\vk{v}'||_q=1, v'_i\geq 0\}$.}

\QED

{\bf Acknowledgement}:Thanks to Enkelejd Hashorva for his suggestions. Thanks to  Swiss National Science Foundation Grant no. 200021-166274.

\bibliographystyle{ieeetr}
\bibliography{Discrete}

\def\lfhook#1{\setbox0=\hbox{#1}{\ooalign{\hidewidth
  \lower1.5ex\hbox{'}\hidewidth\crcr\unhbox0}}}
  \def\polhk#1{\setbox0=\hbox{#1}{\ooalign{\hidewidth
  \lower1.5ex\hbox{`}\hidewidth\crcr\unhbox0}}}
  \def\polhk#1{\setbox0=\hbox{#1}{\ooalign{\hidewidth
  \lower1.5ex\hbox{`}\hidewidth\crcr\unhbox0}}} \def\cprime{$'$}
  \def\cprime{$'$} \def\cprime{$'$}
\begin{thebibliography}{10}

\bibitem{BN1969}
Y.~K. Belyaev and V.~Nosko, ``Characteristics of excursions above a high level
  for a {G}aussian process and its envelope,'' {\em Theory Probab. Appl.},
  vol.~13, pp.~298--302, 1969.

\bibitem{Lindgren1980a}
G.~Lindgren, ``Extreme values and crossing for the chi-square processes and
  other functions of multidimensional {G}aussian process, with reliability
  applications,'' {\em Adv. Appl. Probab.}, vol.~12, pp.~746--774, 1980.

\bibitem{Lindgren1980b}
G.~Lindgren, ``Point processes of exits by bivariate {G}aussian random process
  and extremal theory for the $\chi^2$-processes and its concomitants,'' {\em
  J. Multivariate Anal.}, vol.~10, pp.~181--206, 1980.

\bibitem{Lindgren1989}
G.~Lindgren, ``Slepian model for $\chi^2$-processes with dependent components
  with application to envelope upcrossings,'' {\em J. Appl. Probab.}, vol.~26,
  pp.~36--49, 1989.

\bibitem{Albin1990}
J.~Albin, ``On extremal theory for stationary processes,'' {\em Ann. Probab.},
  vol.~18, no.~1, pp.~92--128, 1990.

\bibitem{Albin1992}
J.~Albin, ``Extremes and crossings for differentiable stationary processes with
  application to {G}aussian processes in $\mathbb{R}^m$ and hilbert space,''
  {\em Stochastic Process. Appl.}, vol.~42, pp.~119--148, 1992.

\bibitem{Berman82}
S.~M. Berman, ``Sojourns and extremes of stationary processes,'' {\em Ann.
  Probab.}, vol.~10, no.~1, pp.~1--46, 1982.

\bibitem{LJ2017}
P.~Liu and L.~Ji, ``Extremes of locally stationary chi-square processes with
  trend,'' {\em Stochastic Processes and their Applications}, vol.~127,
  pp.~497--525, 2017.

\bibitem{ELPNT2017}
L.~Bai, ``Extremes of ${L}^p$-norm of vector-valued {G}aussian processes with
  trend,'' {\em Stochastics}, 2018.

\bibitem{EGCP2018}
L.~Bai, ``Extremes of {G}aussian chaos processes with trend,'' {\em
  https://arxiv.org/abs/1807.00520}, 2018.

\bibitem{KP2018}
K.~I.A. and V.~Piterbarg, ``High excursions of {G}aussian nonstationary
  processes in discrete time,'' {\em Manuscript}, 2018.

\bibitem{Berman92}
S.~M. Berman, {\em Sojourns and {E}xtremes of {S}tochastic {P}rocesses}.
\newblock The Wadsworth \& Brooks/Cole Statistics/Probability Series, Pacific
  Grove, CA: Wadsworth \& Brooks/Cole Advanced Books \& Software, 1992.

\bibitem{MR1342094}
J.~H\"usler, ``A note on extreme values of locally stationary {G}aussian
  processes,'' {\em J. Statist. Plann. Inference}, vol.~45, no.~1-2,
  pp.~203--213, 1995.
\newblock Extreme value theory and applications (Villeneuve d'Ascq, 1992).

\bibitem{MR1062062}
J.~H\"usler, ``Extreme values and high boundary crossings of locally stationary
  {G}aussian processes,'' {\em Ann. Probab.}, vol.~18, no.~3, pp.~1141--1158,
  1990.

\bibitem{MR3679987}
J.~H\"usler and V.~I. Piterbarg, ``On shape of high massive excursions of
  trajectories of {G}aussian homogeneous fields,'' {\em Extremes}, vol.~20,
  no.~3, pp.~691--711, 2017.

\bibitem{ChanLai}
H.~P. Chan and T.~L. Lai, ``Maxima of asymptotically {G}aussian random fields
  and moderate deviation approximations to boundary crossing probabilities of
  sums of random variables with multidimensional indices,'' {\em Ann. Probab.},
  vol.~34, no.~1, pp.~80--121, 2006.

\bibitem{Marek}
M.~Arendarczyk, ``On the asymptotics of supremum distribution for some iterated
  processes,'' {\em Extremes}, vol.~20, no.~2, pp.~451--474, 2017.

\bibitem{MR3638374}
D.~Cheng, ``Excursion probabilities of isotropic and locally isotropic
  {G}aussian random fields on manifolds,'' {\em Extremes}, vol.~20, no.~2,
  pp.~475--487, 2017.

\bibitem{PicandsA}
J.~Pickands, III, ``Upcrossing probabilities for stationary {G}aussian
  processes,'' {\em Trans. Amer. Math. Soc.}, vol.~145, pp.~51--73, 1969.

\bibitem{Pit72}
V.~I. Piterbarg, ``On the paper by {J. Pickands} "upcrosssing probabilities for
  stationary {G}aussian processes",'' {\em Vestnik Moscow Univ Ser. I Mat.
  Mekh. 27, 25-30. English transl. in Moscow Univ. Math. Bull. 1972, 27},
  vol.~27, pp.~25--30, 1972.

\bibitem{debicki2002ruin}
K.~D{\c{e}}bicki, ``Ruin probability for {G}aussian integrated processes,''
  {\em Stochastic Processes and their Applications}, vol.~98, no.~1,
  pp.~151--174, 2002.

\bibitem{DI2005}
A.~Dieker, ``Extremes of {G}aussian processes over an infinite horizon,'' {\em
  Stochastic Process. Appl.}, vol.~115, no.~2, pp.~207--248, 2005.

\bibitem{DE2014}
K.~D\c{e}bicki and K.~Kosi{\'n}ski, ``On the infimum attained by the reflected
  fractional {B}rownian motion,'' {\em Extremes}, vol.~17, no.~3, pp.~431--446,
  2014.

\bibitem{DiekerY}
A.~B. Dieker and B.~Yakir, ``On asymptotic constants in the theory of extremes
  for {G}aussian processes,'' {\em Bernoulli}, vol.~20, no.~3, pp.~1600--1619,
  2014.

\bibitem{DEJ14}
K.~D{\polhk{e}}bicki, E.~Hashorva, and L.~Ji, ``Tail asymptotics of supremum of
  certain {G}aussian processes over threshold dependent random intervals,''
  {\em Extremes}, vol.~17, no.~3, pp.~411--429, 2014.

\bibitem{Pit20}
V.~I. Piterbarg, {\em Twenty Lectures About {G}aussian Processes}.
\newblock London, New York: Atlantic Financial Press, 2015.

\bibitem{Tabis}
K.~D{\c{e}}bicki, E.~Hashorva, L.~Ji, and K.~Tabi{\'s}, ``Extremes of
  vector-valued {G}aussian processes: {E}xact asymptotics,'' {\em Stochastic
  Process. Appl.}, vol.~125, no.~11, pp.~4039--4065, 2015.

\bibitem{DM}
A.~B. Dieker and T.~Mikosch, ``Exact simulation of {B}rown-{R}esnick random
  fields at a finite number of locations,'' {\em Extremes}, vol.~18,
  pp.~301--314, 2015.

\bibitem{SBK}
K.~D\c{e}bicki, S.~Engelke, and E.~Hashorva, ``Generalized {P}ickands constants
  and stationary max-stable processes,'' {\em Extremes}, vol.~20, no.~3,
  pp.~493--517, 2017.

\bibitem{GeneralPit16}
L.~Bai, K.~D\c{e}bicki, E.~Hashorva, and L.~Luo, ``On generalised {P}iterbarg
  constants,'' {\em Methodology and Computing in Applied Probability}, vol.~20,
  pp.~137--164, 2018.

\bibitem{ChCh2018}
L.~Chu and H.~Chen, ``Asymptotic distribution-free change-point detection for
  multivariate and non-euclidean data,'' {\em The Annals of Statistics}, 2018.

\bibitem{Pitchi1994}
V.~Piterbarg, ``High excursions for nonstationary generalized chi-square
  processes,'' {\em Stochastic process. Appl.}, vol.~53, pp.~307--337, 1994.

\bibitem{PL2015}
P.~Liu and L.~Ji, ``Extremes of chi-square processes with trend,'' {\em Probab.
  Math. Statist.}, vol.~36, pp.~1--20, 2016.

\bibitem{Pit96}
V.~I. Piterbarg, {\em Asymptotic methods in the theory of {G}aussian processes
  and fields}, vol.~148 of {\em Translations of Mathematical Monographs}.
\newblock Providence, RI: American Mathematical Society, 1996.

\bibitem{EnkelejdJi2014Chi}
E.~Hashorva and L.~Ji, ``Piterbarg theorems for chi-processes with trend,''
  {\em Extremes}, vol.~18, pp.~37--64, 2015.

\end{thebibliography}
\end{document}